\documentclass[11pt,a4paper]{amsart}

\usepackage[foot]{amsaddr}

\usepackage{amsmath}
\usepackage{amsfonts}
\usepackage{amssymb}
\usepackage{amsthm}
\usepackage{epsfig}
\usepackage{graphicx}
\usepackage{url}
\usepackage[usenames,dvipsnames]{xcolor}
\usepackage[colorlinks=true,linkcolor=Blue,citecolor=Green,urlcolor=black]{hyperref}
\usepackage{nicefrac}
\usepackage{aliascnt}		% differentiation of \autoref names for theorems, lemmas, etc.
\usepackage[alwaysadjust]{paralist}
\usepackage[T1]{fontenc}
\usepackage[utf8]{inputenc}
\usepackage{tikz}
\usetikzlibrary{calc}
\usepackage{enumerate}

\usepackage[numbers]{natbib}

\newtheorem{thm}{Theorem}[section]
\newaliascnt{lem}{thm}
\newaliascnt{cor}{thm}
\newaliascnt{definition}{thm}
\newaliascnt{rem}{thm}
\newaliascnt{prop}{thm}
\newaliascnt{conj}{thm}
\newaliascnt{prob}{thm}

\newtheorem{lem}[lem]{Lemma}
\aliascntresetthe{lem}

\newtheorem*{lem*}{Lemma}

\newtheorem{cor}[cor]{Corollary}
\aliascntresetthe{cor}

\newtheorem*{cor*}{Corollary}

\newtheorem{definition}[definition]{Definition}
\aliascntresetthe{definition}

\newtheorem*{definition*}{Definition}

\newtheorem{rem}[rem]{Remark}
\aliascntresetthe{rem}

\newtheorem*{rem*}{Remark}

\aliascntresetthe{prop}

\newtheorem*{prop*}{Proposition}

\aliascntresetthe{conj}

\newtheorem*{conj*}{Conjecture}

\aliascntresetthe{prob}

\newtheorem*{prob*}{Problem}

\usepackage{dsfont}

\newcommand{\Z}{\mathbb{Z}}

\newcommand{\R}{\mathbb{R}}

\DeclareMathOperator{\poly}{poly}

\DeclareMathOperator{\argmin}{arg.min}

\DeclareMathOperator{\rk}{rk}

\newcommand{\norm}[1]{\|#1\|}

\newcommand{\df}{\mathrel{\mathop:}=}
\newcommand{\fd}{=\mathrel{\mathop:}}

\newcommand{\grad}{\triangledown}

\usepackage[algoruled,vlined]{algorithm2e}
\usepackage{algorithmic}

\numberwithin{equation}{section}

\begin{document}

\title[Convex minimization with hidden low dimension]{Minimizing a low-dimensional convex function over a high-dimensional cube}

\author[C.\ Hunkenschr\"oder]{Christoph Hunkenschr\"oder}
\address[C.\ Hunkenschr\"oder]{Technische Universit\"at Berlin, Institut f\"ur Mathematik, Germany} %Sekretariat MA 4-1, Straße des 17. Juni 136, 10623 Berlin, Germany}
\email{hunkenschroeder@tu-berlin.de}

\date{August 16th, 2022}

\author[S.\ Pokutta]{Sebastian Pokutta}
\address[S.\ Pokutta]{Institute of Mathematics \& AI in Society, Science, and Technology,
Technische Universität Berlin \& Zuse Institute Berlin, Germany}
\email{pokutta@zib.de}

\author[R.\ Weismantel]{Robert Weismantel}
\address[R.\ Weismantel]{ETH Z\"urich, Dept. of Mathematics, Institut f\"ur Operations Research, Switzerland}
\email{weismantel@ifor.math.ethz.ch}

%\thanks{This work was supported by the Swiss National Science Foundation (SNSF) within the project \emph{Convexity, geometry of numbers, and the complexity of integer programming (Nr. 163071)}.}

%\subjclass[2020]{52A41, 90C25, 90C27}

%\keywords{Convex Minimization, Discrete Optimization, Polynomial Time Algorithm}

\begin{abstract}
For a matrix $W \in \Z^{m \times n}$, $m \leq n$, and a convex function $g: \R^m \rightarrow \R$, we are interested in minimizing $f(x) = g(Wx)$ over the set $\{0,1\}^n$.
We will study separable convex functions and sharp convex functions $g$.
Moreover, the matrix $W$ is unknown to us.
Only the number of rows $m \leq n$ and $\|W\|_{\infty}$ is revealed.
The composite function $f(x)$ is presented by a zeroth and first order oracle only.
Our main result is a proximity theorem that ensures that an integral minimum and a continuous minimum for separable convex and sharp convex functions are always ``close'' by.
This will be a key ingredient to develop an algorithm for detecting an integer minimum that achieves a running time of roughly $(m \norm{W}_{\infty})^{\mathcal{O}(m^3)} \cdot \poly(n)$.
In the special case when $(i)$ $W$ is given explicitly and $(ii)$ $g$ is separable convex one can also adapt an algorithm of~\citet{hs}.
The running time of this adapted algorithm matches with the running time of our general algorithm.
\end{abstract}

\thanks{We thank the Hausdorff Research Institute for Mathematics in Bonn for hosting the trimester program \emph{Discrete Optimization}, where large parts of this work developed.\\
This research was partially funded by the DFG Cluster of Excellence MATH+ (EXC-2046/1, project id 390685689) funded by the Deutsche Forschungsgemeinschaft (DFG).\\
The first and the third author were supported by the Einstein Foundation Berlin.}

\maketitle

%\tableofcontents

\section{Introduction}
\label{sec:intro}
Minimizing a convex function over a convex set is a versatile and important problem in many modern research fields, such as discrete and combinatorial optimization, or machine learning, to name a few.

In certain regression problems, the objective is of the form $f(x) = g(Wx)$ for a matrix $W \in \Z^{m \times n}$ with $m \ll n$,
and a low-dimensional convex function $g$.
While this kind of objective can be represented in two distinct parts ($W$ and $g$),
we assume that we only have oracle access to the composition $f(x)$.
%that are relatively easy to handle for themselves, we still might have only oracle access to the composition $f(x)$.

The focus of this paper is on solving
\[
\min \{g(Wx) : \ x \in \{0,1\}^n\},
\]
for some integer matrix $W \in \Z^{m \times n}$, and convex function $g: \R^m \rightarrow \R$.
We emphasize that the matrix $W$ is unknown to us.
Only (bounds on) the number of rows $m \leq n$ and $\norm{W}_{\infty}$ are revealed.
When it comes to $g$,
we assume that the gradient $\grad g(u)$ is integral whenever $u$ is, and that 
$g$ has Lipschitz continuous gradients (Definition~\ref{def:smooth}).
This can be viewed as a normalization for the oracle presentation of $g$.
In addition to these general assumptions, we highlight two classes of convex functions.
We will assume that $g$ is either separable convex (Section~\ref{sec:proximity-sep-conv}) or sharp (Section~\ref{sec:proximity-sharp}).
In both cases, we will establish that an integral minimum can be found near a continuous minimum; such a \emph{proximity result} is crucial for our algorithm.
If a similar result can be shown for another class of convex functions with Lipschitz continuous gradients, our results readily extend to this class.

An example for our setting stems from \emph{compressed sensing} where the objective $f(x) = g(Wx) = \| b - Wx\|^2$ is considered, cf.~\citet{braun2015info}.
The vector $x \in \R^n$ is a high-dimensional hidden signal we want to reconstruct, while the vector $b \in \Z^m$ is the outcome of several measurements that were conducted.
Specifically, $m$ is the number of measurements and therefore known to the observer, while $W$ is not explicitly given, but only implicitly through a zeroth and first order oracle of $f$.
%Depending on the specific setting, the matrix $W \in \Z^{m \times n}$ contains data related to the measurements, for instance discrete energies.
If $W$ and $x$ are integral, so is $\grad f(x) = -(b - Wx)^\intercal W$.
Moreover, we can estimate $\|\grad g(u) - \grad g(v) \|_\infty = 2 \|u-v\|_{\infty}$, hence the Lipschitz constant of the gradient is bounded by $2$.

We develop an algorithm that does not require access to $W$ and $g$ separately, but only requires zeroth and first order access to $f(x) = g(Wx)$, i.e., it queries the value and the gradient of $f$ for a given point $x \in [0,1]^n$.
Informally, our main result reads as follows.
\begin{thm}[informal, cf.\ Cor.~\ref{cor:sep-conv} and~\ref{cor:sharp-conv}]
\label{thm:informal}
Let $f(x) = g(Wx)$ for $W \in \Z^{m \times n}$, and $g$ a convex function with Lipschitz continuous gradients that is either separable or sharp.

Given bounds $\bar{m} \geq m$ and $\Delta \geq \norm{W}_\infty$, and a subroutine that computes for a face $F \subseteq [0,1]^n$ a continuous constrained optimum for $\min \{f(x): x \in F\}$,
we can find an integer optimum $z^\star \in \{0,1\}^n$ 
in
\[
(\bar{m} \Delta + 1 )^{\mathcal{O}(\bar{m} m^2)} \cdot poly(n)
\]
arithmetic operations and calls to a zeroth and first order oracle of $f$, plus $2n$ calls to the subroutine.
\end{thm}

\subsection{Related literature}
\label{sec:literature}
For continuous optimization, it is a common approach to use the gradients for minimizing a convex objective $f$.
We refer the reader to the surveys of~\citet{cg_survey} and~\citet{Nocedal}.

In integer programming, it is less obvious to design an iterative gradient based algorithm.
Lenstra's algorithm can be used to check whether the convex set $C_{\beta} = \{x \in [0,1]^n: \ f(x) \leq \beta \}$ contains an integer point.
Then, a binary search on $\beta$ eventually leads to the minimum (cf.~\cite[Theorem 6.7.10]{gls}).
While this approach readily allows for additional (convex) constraints on the set, its running time is exponential in the dimension $n$.
%, and certain assumptions on the encoding of $g$ need to be made.

\citet{oertel2014integer} proposed an approach using the gradients of $f$.
They start with considering a polyhedron that must contain the minimal solution.
Then they find an integer point deep inside the polyhedron, and use its gradient as a cutting plane.
Reducing the volume of the polyhedron by a constant factor, they iterate until they are left with a polytope of small volume.
As a consequence, there is an integral vector $a$ so that only a few of the hyperplanes $H_z = \{x: \ a^\intercal x = z\}$, $z \in \Z$ intersect the polytope non-trivially, and they can recurse in lower dimension.
If the dimension is fixed, they reduce the problem in this way to solving a polynomial number of mixed integer linear programs.

A quite beautiful approach in two dimensions was presented by~\citet{paat2021constructing}.
They start with a \emph{gradient polyhedron} $P$, i.e., $P$ is described by the inequalities $\grad f (z_i)^\intercal x \leq \grad f (z_i)^\intercal z_i$ induced by some integer points $z_1,\dots,z_m \in \Z^2$.
If this polyhedron does not contain any integer point in its interior, then one of $z_1,\dots,z_m$ has to be optimal.
Otherwise, if there exists an integer point $z \in \text{interior} (P)$, they iteratively find it and swap it against one of $z_1,\dots,z_m$, yielding a new gradient polyhedron $P^\prime$.
A discrete analogue of Helly's theorem, proven by~\citet{doignon1973convexity}, implies that $m=4$ suffices.

\citet{de2009convex} are interested in \emph{maximizing} a convex function $f(x) = g(Wx)$, $W \in \Z^{m \times n}$, over a system $Ax = b$, $x \in \Z^n_{\geq 0}$.
The authors use the Graver basis of $A$ to derive a polynomial-time algorithm, if $m$ and $\norm{A}_{\infty}$ are constant.
An earlier work, considering a more abstract setting, is by~\citet{onn2004convex}.
Optimizing $g(Wx)$ over a polyhedron $P \subseteq \R^n$ for an \emph{arbitrary} function $g$ was also considered by~\citet{adjiashvili2015polyhedral}, though they focus on the structural relation between the projections $W(\Z^n \cap P)$ and $WP \cap \Z^m$.
In particular, they already assume an optimization oracle for $g$.
These works are, however, the only ones considering the concatenation $f(x) = g(Wx)$ of a projection and a low-dimensional convex function, we are aware of.

Closest to our paper is an approach of~\citet{hs} that applies to separable convex functions provided that $W$ is given explicitly.
Their main results include a proximity theorem and an algorithm to compute the integral minimum.
\citet{hs} extend a result of~\citet{sensitivity} for linear functions to separable convex functions $g$.
They show that if $x$ is a fractional minimal solution to $\min \{g(x) : \ Ax \leq b\}$ with $A \in \Z^{k \times n}$, then there exists a minimal solution $z^\star$ to the corresponding integer nonlinear program satisfying $\norm{x^\star - z^\star}_{\infty} \leq n \Delta$, where $\Delta$ is the largest subdeterminant of $A$.
This implies $\delta = \| x^\star - z^\star \|_1 \leq n^2 \Delta$.
 
The minimum $z^\star$ has to be in a cube around the origin of edge length $2 B$ for some $B$ depending on the input data.
% = 2 n \norm{W}_\infty$.
\citet{hs} develop a framework that reduces the problem to minimizing a sequence of linear (integer) programs.
In a nutshell, they replace each $g_i$ by a piecewise linear function $\tilde{g_i}$ on the box $[-B,B]^m$, with only a ``few'' pieces.
They split each $z_i$ into several variables, one per piece of $\tilde{g_i}$, yielding a linear program in a higher dimension.
As each $\tilde{g_i}$ only has a small number of pieces, the solution $z^\star$ to the modified problem might not be minimal for the initial problem,
but they show that the minimal solution to the initial problem is within a box $z^\star + [-B/2,B/2]^m$, allowing them to iterate this process.
After $\log(B)$ steps, they are left with a fairly small box, on which they solve an integer program.
In the original work, the number of pieces required for $\tilde{g_i}$, as well as the achieved running times, are presented in dependence on the largest subdeterminant of $A$ and $n$.

In order to adapt their approach to our setting one needs some modifications that we outline next.
% into an algorithm for our problem consider the reformulation
Given that $W$ is known explicitly, we may consider the following reformulation:
\begin{align*}
%\label{eq:reformulation}
\min \{\sum_{i=1}^m &g_i(y_i): \ Wx = y, \ x \in \{0,1\}^n \}.
\end{align*}
It is not difficult to verify that we can assume without loss of generality that all columns of $W$ are distinct (by possibly changing the variable bounds $[0,1]$).
The number of distinct columns is upper bounded by $(2 \|W\|_{\infty} + 1)^m$ and hence the dimension $n$ is at most $(2 \|W\|_{\infty} + 1)^m + m$.
Moreover, using Hadamard's inequality one can upper bound the largest subdeterminant of $W$ by $\Delta \leq (m\|W\|_{\infty})^m$.
Hence, with these adaptations the proximity theorem of~\citet{hs} gives
\begin{align}
\label{eq:hs-proximity}
\| x^\star - z^\star \|_1 \leq n^2 \Delta &\leq (2m \|W\|_{\infty} + 1)^{3m}.
\end{align}
In Section~\ref{sec:proximity-sep-conv} we will improve the exponent in~\eqref{eq:hs-proximity}.

In order to adapt their algorithm to this setting one needs to apply an integer programming algorithm for linear objectives in standard form representation.
With the current improvements on dynamic programming by~\citet{ew} one can analyze their running time and obtain 
$(m \norm{W}_\infty)^{\mathcal{O} (m^2)} \poly(n)$.

Alternatively, the above reformulation can also be solved with Graver basis augmentation techniques (see e.g.,~\cite{de2012algebraic, hemmecke2011polynomial, onn2010nonlinear}), yielding similar running times.
However, we emphasize that for the initial reformulation, explicit knowledge of $W$ is necessary, which is in contrast to the algorithm that we will develop later.
If $W$ is given to us, our algorithm can be sped up to the same asymptotic running time, cf.\ Remark~\ref{rem:exponent}.

\subsection{Example in one dimension: solving the knapsack problem}
We illustrate the idea for our approach applied to the knapsack problem:
Let $W \in \Z^{1 \times n}$, and $T \in \Z_{> 0}$.
Does there exist $x \in \{0,1\}^n$ with $Wx = T$?

Defining $g(y) = (T-y)^2$, $f(x) = g(Wx)$,
the original instance is a yes-instance, if and only if the minimum value is zero.

We generate a sequence $(x^j)_{j=0}^N \subseteq \{0,1\}^n$ of points as follows, where $x^0 \in \{0,1\}^n$ is an arbitrary starting point.
If $(\grad f(x^j))_i \geq 0$ whenever $x^j_i = 0$, and $\grad f(x^j)_i \leq 0$ whenever $x^j_i = 1$, then the condition
\[
f(x) \geq f(x^j) + \grad f(x^j)^\intercal (x-x^j)
\]
for each $x \in \{0,1\}^n$ certifies optimality of $x^j$.

Otherwise, let $k$ be an index violating this condition.
We obtain $x^{j+1}$ from $x^j$ by applying the map $x_k \mapsto 1 - x_k$, leaving the other coordinates untouched.
Recall that $W \in \Z^{1 \times n}$, implying $\grad g (Wx) \in \R$.
Consequently, writing $\grad f(x)^\intercal = \grad g(Wx) W$ shows that the gradients of $f$ at $x^j$ and $x^{j+1}$ are linearly dependent.
Hence, we can write $\grad f(x^{j+1}) = \alpha_j \grad f(x^j)$ for some $\alpha_j \in \R$.
Depending on the sign of $\alpha_j$, we distinguish three cases.

If $\alpha_j = 0$, then $\grad f(x^{j+1}) = 0$ and $x^{j+1}$ is optimal.

If $\alpha_j > 0$, either $x^{j+1}$ is minimal or we obtain a point $x^{j+2}$, following the argumentation above.
Observe that $\grad f(x^j)$ and $\grad f(x^{j+1})$ are sign compatible, implying that no index is chosen more than once, and after at most $n$ steps, $x^{j+1}$ is either minimal, or we obtain $\alpha_j < 0$.

If $\alpha_j < 0$, the gradient changed its sign, i.e., we walked too far.
Any optimal point $z^\star \in \{0,1\}^n$ satisfies 
\[
\grad f(x^j)^\intercal z^\star \leq \grad f(x^j)^\intercal x^j,
\quad \text{and} \quad
\grad f(x^{j+1})^\intercal z^\star \leq \grad f(x^{j+1})^\intercal x^{j+1}
\]
\[
\Leftrightarrow \grad f(x^j)^\intercal x^{j+1} \leq \grad f(x^j)^\intercal z^\star \leq \grad f(x^j)^\intercal x^j.
\]
In other words, the gradients of $f$ at $x^j$ and $x^{j+1}$ describe a small polyhedron containing a minimal solution.

Let us assume the gradients are integral and bounded, for instance, $\norm{\grad f(x^j)}_\infty$, $\norm{\grad f(x^{j+1})}_\infty \leq 2 \norm{W}_\infty^2$.
Since $\norm{x^{j+1} - x^j}_1 = 1$, the polyhedron is intersected by at most $2 \norm{W}_\infty^2$ hyperplanes on which $z^\star$ may lie.

In total, we construct a sequence $(x^j)_{j=1}^N$ of at most $N \leq n$ points, and are left with at most $2 \norm{W}_\infty^2$ hyperplanes for which we have to check whether they contain a point in $\{0,1\}^n$.
If a hyperplane $H$ contains several points $z^i,z^j \in H$, observe that $f(z^i) = f(z^j)$ for all these points.
Therefore, we can iterate over all hyperplanes and pick a minimal point.
This is done with an integer programming algorithm we will be using as a blackbox.

As a final remark, this example already shows that enumerating hyperplanes (or, more generally, affine subspaces), can (most likely) not be avoided.
If there was an algorithm polynomial in $n$ and not depending on $\norm{W}_{\infty}$, we could solve the knapsack problem in polynomial time.

\subsection{The approach in higher dimension}
\label{sec:contributions}
When $m \geq 2$, we lose the notion of the gradient changing its sign, as different gradients can be linearly independent.
This makes it more difficult to define improving steps which converge to a small set around the minimum.
Instead, the gradients will help restricting the search space in another way.

Informally, we proceed as follows.
The algorithm takes as input a fractional (constrained) minimal solution $x^\star$ for $f(x) = g(Wx)$, $x\in [0,1]^n$.
Our first guess for an integer minimum $z^\star$ will be $\lceil x^\star \rfloor$, where we round component-wise.
During the course of the algorithm, we will construct a set $S \subseteq \{0,1\}^n$; initially, it will only contain the point $\lceil x^\star \rfloor$.
The points $z \in S$ will provide linearly independent gradients $h_z \df \grad f(z)$ such that $z^\star$ will satisfy $h_z^\intercal z^\star \leq h_z^\intercal z$ for all $z \in S$.% for some $\beta_i \in \Z_{\geq 0}$.

Since the gradients are integral,
we can iterate over the values of the dual gaps, 
i.e., the affine subspaces $\{x \in \R^n: \ h_z^\intercal x = h_z^\intercal z - \beta_z \ \forall z \in S\}$,
$\beta_z \in \Z_{\geq 0}$ and check whether they contain a new integer point $z^\prime \in \{0,1\}^n$.
It turns out that we can either add $z^\prime$ to the set $S$, effectively reducing the dimension of the affine subspaces, or $z^\prime$ is a candidate for a minimum.

How do we limit the number of considered subspaces though?
The key for this are \emph{proximity results}: We show that if the smooth convex function $g$ is either separable convex or sharp, then $z^\star$ must be close to $x^\star$.
However, even if they are close, a large gradient $h_z$ could still imply that we have to consider a vast number of hyperplanes, as $h_z^\intercal x^\star$ and $h_z^\intercal z^\star$ can differ substantially.
To circumvent this, we show that if $\grad f(x^\star)$ is large, certain variables of $z^\star$ are already forced to either $0$ or $1$, leaving us with a lower-dimensional problem where all gradients of interest will have small entries.

The remainder of the paper is structured as follows.
In Section~\ref{sec:prelims}, we introduce necessary concepts and prove several structural results.
A key notion is $\delta$-proximity for a function $f$, presented in Definition~\ref{def:proximity}.
In Section~\ref{sec:algo}, we present our algorithm in pseudo-code and analyze its performance in dependency on the parameter $\delta$.
In Section~\ref{sec:proximity}, we will analyze proximity questions for two different classes of convex functions, leading to concrete upper bounds on $\delta$ and explicit running times in these cases.

\section{Preliminaries}
\label{sec:prelims}
Throughout this paper, we keep the assumptions made in the introduction, i.e., $W \in \Z^{m \times n}$, and for $g: \R^m \rightarrow \R$ convex, 
$\grad g (u) \in \Z^m$ whenever $u \in \Z^m$.
%$\grad g(u) \in \tfrac{1}{c(u)}\Z^m$ for a constant $c(u)$ whenever $u \in \Z^m$.
We also require $g$ to have Lipschitz continuous gradients on $[0,1]^n$, which we specify now.

In the literature, $g$ having Lipschitz continuous gradients usually means that there is a constant $L > 0$ such that 
$\norm{\grad g (x) - \grad g (y)}_2 \leq \norm{x-y}_2 L$ holds for all $x,y \in \R^m$.
However, in the field of integer programming, it is more natural to work with the $\ell_1$- and $\ell_{\infty}$-norms.
Clearly, all these norms are related by factors of at most $m$.
As other estimates will dominate this factor, we settle for a definition in terms of the $\ell_1$- and $\ell_{\infty}$-norms.
\begin{definition}[Lipschitz continuous gradients]
\label{def:smooth}
Let $g: \R^m \rightarrow \R$ be a continuously differentiable convex function, and $L > 0$.
We say that $g$ has $L$-Lipschitz continuous gradients, if 
\begin{align}
\label{eq:smooth}
\norm{\grad g (u) - \grad g(v)}_\infty &\leq L \norm{u-v}_1 \qquad \forall \ u,v \in \R^m.
\end{align}
\end{definition}
Note that a scaling of $g$ affects the Lipschitz constant $L$.
In order to make the representation of $g$ and $L$ unique, we assume integral gradients $\grad g(x)$ for $x \in \Z^n$.
This can be viewed as a normalization of a scaling of $g$ that seems natural when evaluating $\grad g$ at integral points.
A standard result for convex minimization that we will use frequently is the following lemma, together with a simple implication.
\begin{lem}[First order optimality condition]
\label{lem:first-order-opt}
Let $g: \R^m \rightarrow \R$ be convex and differentiable, and $v^\star \in P$ minimal for $g$ over a convex set $P \subseteq \R^m$.
For all $u \in P$, we have $\grad g (v^\star)^\intercal (v^\star - u) \leq 0$.
\end{lem}
\begin{cor}
\label{cor:first-order-opt}
Let $f: [0,1]^n \rightarrow \R$ be convex and differentiable, and $x^\star \in [0,1]^n$ a continuous minimal solution for $f$.
The following two implications hold.
\begin{align*}
\grad (g(x^\star))_k < 0 &\Rightarrow x^\star_k = 1 \\
\grad (g(x^\star))_k > 0 &\Rightarrow x^\star_k = 0 \\
\end{align*}
\end{cor}

As mentioned before, a crucial argument for our algorithm to work is that an integral constrained minimum of $f$ is close to a fractional constrained minimum of $f$.
In general, we cannot hope for a bound on the $\ell_1$-norm independent on the dimension; the unique continuous minimizer of $f(x) = \sum_{i=1}^n (2x_i - 1)^2$ is $\tfrac{1}{2} \cdot \mathbf{1}$ which has distance $\mathcal{O}(n)$ to $\Z^n$.

A simple observation shows that in our setting however, there is always a minimum with only a few fractional entries, ruling out this example.
\begin{lem}[Few fractional entries]
\label{lem:frac-entries}
Let $f(x) = g(Wx)$ with $W \in \Z^{m \times n}$, and $\hat{x} \in [0,1]^n$.
There exists $x^\prime \in [0,1]^n$ with at most $m$ fractional entries, such that $f(x^\prime) = f(\hat{x})$.
\end{lem}
\begin{proof}
Let $\hat{x} \in [0,1]^n$, and let $x^\prime$ be a vertex of the polytope $P \df \{x \in [0,1]^n : \  W x = W \hat{x} \}$.
Clearly, we have
$
f(x^\prime) = g(Wx^\prime) = g(W\hat{x}) = f(\hat{x}).
$
By basic LP theory, the number of rows of $W$ is an upper bound on the number of fractional entries in $x^\prime$.
This shows the claim.
\end{proof}
The proof reveals that given any fractional minimum, we can even compute a minimum with at most $m$ fractional entries, provided that we know $W$.
If we do not know $W$ we will have to assume that there is an algorithm at our disposal which solves the underlying continuous minimization problem.
This will be discussed in Lemma~\ref{lem:frac-opt}, but for the benefit of presentation, we assume such a minimum to be given.

Bearing this in mind, we can formalize the proximity of a constrained fractional and a constrained integral minimal solution.
\begin{definition}[$\delta$-proximity]
\label{def:proximity}
Let $f(x) = g(Wx)$ with $g: \R^m \rightarrow \R$ convex and $W \in \Z^{m \times n}$.
We say that $f$ \emph{admits a $\delta$-proximity}, if there is a constant $\delta \in \R_{\geq 0}$ such that the following holds:
For every $x^\star \in \argmin \{f(x): \ x \in [0,1]^n \}$ with at most $m$ fractional entries, there exists $z^\star \in \argmin \{f(z): \ z \in \{0,1\}^n\}$ with $\norm{x^\star - z^\star}_1 \leq \delta$.
\end{definition}

We will eventually show that for specific convex functions $g$, the function $f(x) = g(Wx)$ admits a $\delta$-proximity for small values of $\delta$ depending on $W$, but independent of $n$.
As the proofs differ for separable convex and sharp convex functions, we defer to Theorems~\ref{thm:proximity-sep-conv} and~\ref{thm:proximity-sharp}.

In the course of the algorithm, we will encounter the following setting.
A solution $z^\prime$ to a system $Hz = b^\prime$ is given, and we are interested in a solution $z^{\prime \prime}$ for a slightly altered right-hand side, $Hz = b^{\prime \prime}$.
A result of~\citet{sensitivity} shows that if $b^{\prime \prime}$ is close to $b^\prime$ and the latter system has an integer solution, then we can find one, $z^{\prime \prime}$ say, that is close to $z^\prime$.
In the field of integer programming, relating $\norm{z^\prime - z^{\prime \prime}}$ to $\norm{b^\prime - b^{\prime \prime}}$ is usually referred to as a \emph{sensitivity result}.
While~\citet{sensitivity} study polyhedra in inequality representation, we provide a sensitivity result for systems in standard form, which is appropriate for our setting.
In contrast to~\citet{sensitivity}, we use techniques inspired by~\citet{ew} and obtain the following.

\begin{lem}[Sensitivity for integer programming]
\label{lem:sensitivity}
Let $H \in \Z^{k \times n}$, and $b^\star \in \R^k$, $b \in \Z^k$ such that there exists $x^\star \in [0,1]^n$ with $Hx^\star = b^\star$.
By Lemma~\ref{lem:frac-entries}, assume that $x^\star$ has at most $\varphi \leq m$ fractional entries.
If the system $Hz = b$, $z \in \{0,1\}^n$ has a solution, then there exists a solution $\hat{z} \in \{0,1\}^n$ with 
\[
\norm{x^\star - \hat{z}}_1 \leq \left( \left\lceil \tfrac{\norm{b^\star - b}_\infty}{\norm{H}_\infty}  \right\rceil + \varphi + 1 \right) (2k \norm{H}_\infty + 1)^k.
\]
\end{lem}
\begin{proof}
For ease of notation, we will write $\Delta \df \norm{H}_\infty$.
Let $\lceil b^\star \rfloor_b$ be $b^\star$ rounded towards $b$, i.e., $(\lceil b^\star \rfloor_b)_i = \lfloor b^\star_i \rfloor$ if $b_i \leq b_i^\star$, and 
$(\lceil b^\star \rfloor_b)_i = \lceil b^\star_i \rceil$ otherwise.
Let $\{ b^\star \}_b = b^\star - \lceil b^\star \rfloor_b$, and define $\lceil x^\star \rfloor_z$ and $\{x^\star\}_z$ similarly, where $z$ is a solution to $Hz = b$, $z \in \{0,1\}^n$ that is closest to $x^\star$ (in terms of $\|\cdot\|_1$).
Consider 
\begin{align*}
H(x^\star - z) &= b^\star - b \\
\Leftrightarrow \quad \underbrace{ H(\lceil x^\star \rfloor_z - z)}_{\fd u} + \underbrace{H \{x^\star\}_z - \{b^\star\}_b}_{\fd v} + \underbrace{b - \lceil b^\star \rfloor_b}_{\fd w} &= 0, 
\end{align*}
and observe that since $u,w \in \Z^k$ we also have $v \in \Z^k$.
Since $x^\star$ has at most $\varphi$ fractional entries, we have $\norm{v}_\infty \leq \varphi \Delta$.
Therefore, we can write $v = v_1 + \dots + v_{\varphi}$ for some vectors $v_i \in \Z^k$ with $\norm{v_i}_\infty \leq \Delta$.
Observe that we can write $w = w_1 + \dots + w_{\ell}$ for some vectors $w_i \in \Z^k$ such that $\norm{w_i}_\infty \leq \Delta$ and $\ell \leq \left\lceil \tfrac{\norm{b^\star - b}_\infty}{\Delta} \right\rceil$.

We next construct a number $N$ and a decomposition $u = u_1 + \dots + u_N$ in the following manner:
If $(\lceil x^\star \rfloor_z - z)_i = 1$, we add the $i$-th column $H_i$ of $H$ to the sum, if $(\lceil x^\star \rfloor_z - z)_i = -1$, we add $-H_i$ of $H$ to the sum, and if $(\lceil x^\star \rfloor_z - z)_i = 0$, we do not add a vector to the sum.
Hence, $N = \norm{\lceil x^\star \rfloor_z - z}_1$, and $\norm{u_i}_\infty \leq \Delta$ for each $i$.

Combining the three decompositions of $u,v$, and $w$, we obtain $N + \varphi + \ell$ vectors satisfying $\norm{u_i}_\infty,\norm{v_i}_\infty,\norm{w_i}_\infty \leq \Delta$, and summing up to zero,
\begin{align}
\label{eq:sum}
\sum_{i=1}^N u_i + \sum_{i=1}^{\varphi} v_i + \sum_{i=1}^{\ell} w_i &= 0.
\end{align}
We will now apply a Lemma attributed to~\citet{steinitz1913bedingt} and improved by~\citet{grinberg1980value} to this sum, stating the following:
Given $\sum_{i=1}^M x_i = 0$ with $x_i \in \R^k$ and $\norm{x_i} \leq 1$ for an arbitrary norm $\| \cdot \|$, there is a permutation $\pi: \{1,\dots,M\} \rightarrow \{1,\dots,M\}$ such that $\norm{\sum_{i=1}^{r} x_{\pi(i)}} \leq k$ for each $r \in [M]$.
As the statement is invariant under scaling, we can reorder (and rename) the sum~\eqref{eq:sum} to 
\[
\sum_{i=1}^{N + \varphi + \ell} a_i = 0 \qquad \text{so that} \qquad
\left\| \sum_{i=1}^{r} a_i \right\|_\infty \leq k \Delta
\]
for each $r \leq N + \varphi + \ell$.
Thus, there are at most $(2k \Delta + 1)^k$ possible values for any prefix sum.
In particular, if $N + \varphi + \ell > (2 k \Delta + 1)^k$, there exist two indices $j_1 < j_2$ with $j_2 - j_1 \leq (2k \Delta + 1)^k$, such that
\[
\sum_{i = 1}^{j_1} a_i = \sum_{i=1}^{j_2} a_i
\quad \Longleftrightarrow \quad
\sum_{i=j_1+1}^{j_2} a_i = 0.
\]
The latter expression, i.e., a set of elements summing up to $0$ is called a \emph{cycle}.
Even more, if $N + \varphi + \ell > (\varphi + \ell + 1) (2k \Delta + 1)^k$, there are at least $\varphi + \ell + 1$ cycles.
This can be seen by choosing a shortest cycle (ensuring $j_2 - j_1 \leq (2 k \Delta + 1)^k$), deleting it, and searching for a cycle in the remaining sequence.
As there are only $\varphi + \ell$ elements $a_i$ corresponding to a vector $v_{\ell}$ or $w_{\ell}$, 
there is a cycle $S$ such all vectors $a_k$, $k \in S$ correspond to elements in $\{ u_1,\dots,u_N \}$.

Note that if $u_{i} = - H_{\ell}$ for some $\ell$, this implies $(\lceil x^\star \rfloor_z - z)_{\ell} = -1$, i.e., $x^\star_{\ell} = 0$ and $z_{\ell} = 1$.
Hence, deleting $u_i$ corresponds to switching $z_{\ell}$ to $0$.
Similarly, if $u_i = H_j$, deletion corresponds to switching $z_j$ to $1$.
In both cases, $\norm{x^\star - z}_1$ decreases.
It follows that if we delete a cycle that consists only of elements in $\{u_1,\dots,u_N\}$, we obtain a new vector $z^\prime$ such that $H z = H z^\prime$, and $\norm{x^\star - z^\prime}_1 < \norm{x^\star - z}_1$.
This contradicts the choice of $z$.
Recalling $N = \norm{\lceil x^\star \rfloor_z - z}_1$ we therefore have
\[
\norm{x^\star - z}_1 \leq N + \varphi \leq (\varphi + \ell + 1) (2k \Delta + 1)^k - \ell. \qedhere
\]
\end{proof}
\begin{rem}
From Lemma~\ref{lem:sensitivity}, we obtain a sensitivity result for two integer solutions by setting $\varphi = 0$.
\end{rem}

Recall Corollary~\ref{cor:first-order-opt}: If $x^\star \in [0,1]^n$ is a continuous minimum, knowing $\grad f(x^\star)$ allows us already to fix all entries $x^\star_k$ to either $0$ or $1$ whenever $(\grad f(x^\star))_k \neq 0$.
It appears natural to ask if there are similar conditions on $\grad f(x^\star)$ that allow us to reconstruct entries in an integer minimum $z^\star$ as well.
As it turns out, we can strengthen this condition:

\begin{lem}[Small gradients]
\label{lem:small-grads}
Let $f(x) = g(Wx)$ admit a $\delta$-proximity for some $\delta \in \R_{\geq 0}$, and let $g$ be convex and have $L$-Lipschitz continuous gradients.
For any constrained continuous minimum $x^\star \in [0,1]^n$ of $f$ with at most $m$ fractional entries, there exists an integer minimum $z^\star \in \{0,1\}^n$ with the following property.
\begin{enumerate}
\item If $(\grad f(x^\star))_k > \max \{(\delta - 1) L, 0\}$, then $z^\star_k = 0$.
\item If $(\grad f(x^\star))_k < - \max\{(\delta - 1) L, 0 \}$, then $z^\star_k = 1$.
\end{enumerate} 
\end{lem}
\begin{proof}
As the two claims are symmetric, we will only show Claim $(1)$.
For a constrained continuous minimum $x^\star \in [0,1]^n$ with at most $m$ fractional entries, let $k$ be an index with 
$(\grad f (x^\star))_k > \max \{(\delta - 1) L, 0\}$.
By Corollary~\ref{cor:first-order-opt} we have $x^\star_k = 0$.
Let $z^\star$ be an integer minimum with $\norm{x^\star - z^\star}_1 \leq \delta$ by proximity.
If $\delta < 1$, we have $|x^\star_k - z_k^\star| \leq \| x^\star - z^\star \|_1 < 1$, hence we can conclude $z^\star_k = x^\star_k = 0$.

Otherwise assume for the sake of contradiction that $z^\star_k = 1$, which implies that for $\hat{z} \df z^\star - e_k$ we have $\norm{x^\star - \hat{z}}_1 \leq \delta - 1$.
Moreover, we have $\norm{\grad f (x^\star) - \grad f (\hat{z})}_{\infty} \leq \norm{x^\star - \hat{z}}_1 L$ and therefore
\[
(\grad f (\hat{z}))_k \geq (\grad f (x^\star))_k - (\delta - 1) L > (\delta - 1) L - (\delta - 1) L = 0.
\]
Thus, $f(z^\star) \geq f(\hat{z}) + ( \grad f(\hat{z}) )^\intercal (z^\star - \hat{z}) > f(\hat{z})$,
which is a contradiction.
\end{proof}
\begin{rem}
\label{rem:prepro}
Lemma~\ref{lem:small-grads} allows us to perform some preprocessing steps.
In fact, we can set specific variables to either $0$ or $1$ if the assumptions of Lemma~\ref{lem:small-grads} apply.

To see this, let $x^\star \in [0,1]^n$ be a constrained fractional minimal solution for $f$ with at most $m$ fractional entries.
Assume w.l.o.g.\ that the first $n-k$ entries of $\grad f (x^\star)$ are small in absolute value (i.e., at most $(\delta -1)L$), and the last $k$ entries are large.
We can split $W = (W_1, W_2)$ and $x^\star = \binom{x^\star_1}{x^\star_2}$ accordingly.
By Corollary~\ref{cor:first-order-opt}, $x^\star_2 \in \{0,1\}^k$.
We can define $g^\prime(v) = g(v + W_2 x^\star_2)$, which is a translate of $g$.
As our assumptions on $g$ are invariant under translation, we can minimize $f^\prime (x^\prime) \df g^\prime(W_1x^\prime)$ instead of $f$.
A minimal solution $z^\star$ to the modified problem can be extended to a minimal solution
$\binom{z^\star}{x_2^\star}$ of the initial problem eventually, using Lemma~\ref{lem:small-grads}.
However, the gradient $\grad f^\prime(x^\star_1)$ is simply the gradient $\grad f(x^\star)$ restricted to the first $n-k$ entries, hence $\norm{\grad f^\prime(x^\star_1)}_\infty \leq (\delta - 1)L$.
\end{rem}

\section{The algorithm}
\label{sec:algo}
Now that we have all necessary concepts at our disposal, let us state and analyze the algorithm.
We assume that $W \in \Z^{m \times n}$, the function $g$ has $L$-Lipschitz continuous gradients that satisfy $\grad g(u) \in \Z^m$ whenever $u \in \Z^m$, and $f(x) = g(Wx)$ admits a $\delta$-proximity for some $\delta \geq 0$.
The algorithm has access to $f$ via a zeroth and first order oracle.
As input, we assume %$m$, $\norm{W}_{\infty}$, and
a fractional minimum $x^\star \in [0,1]^n$ of $f$ with at most $m$ fractional entries and the constants $\delta,L$.
The output will be $z^\star \in \{0,1\}^n$ minimizing $f$.

\bigskip
\begin{algorithm}[H]
\begin{algorithmic}[1]
 \REQUIRE {A number $\delta \in \R_{\geq 0}$ s.t.\ $f$ admits a $\delta$-proximity, a Lipschitz constant $L$ of $\grad g$, and a continuous constrained minimum $x^\star \in [0,1]^n$ with at most $m$ fractional entries.}
 \ENSURE {An integer minimum $z^\star \in \{0,1\}^n$.}

 \IF {$\norm{\grad f(x^\star)}_{\infty} > \max\{ (\delta-1)L, 0\}$}
 	\STATE Apply the preprocessing discussed in Remark~\ref{rem:prepro}. \label{alg:restrict-to-face}
 \ENDIF
 \STATE $\hat{z} \df \lfloor x^\star \rceil$
 \IF {$\grad f (\hat{z}) = 0$}
 	\RETURN {$(\hat{z})$}
 \ENDIF
 \STATE $\Gamma \df \{ \grad f (\hat{z}) \}$  % ^\intercal \in \Z^{1 \times n}$  
 \REPEAT  \label{alg:outerloop}
 \STATE $r \leftarrow |\Gamma|$
 \STATE Let $H \in \Z^{r \times n}$ be the matrix whose rows are the vectors in $\Gamma$. \label{alg:end-init}
 \STATE $b^\star \leftarrow H x^\star \in \R^{r}$
 \FORALL{$b \in \{x \in \Z^{r}: \norm{x - b^\star}_\infty \leq \norm{H}_\infty \delta \}$} \label{alg:b-loop}
	\IF{$\min \{ \norm{x^\star - z}_1: Hz = b, z \in \{0,1\}^n\}$ has a solution $z_b$ with $\norm{z_b - x^\star}_1 \leq \delta$} \label{alg:ip} 

			\IF{$\grad f(z_b)$ is linearly independent of $\Gamma$}
			\STATE $\Gamma = \Gamma \cup \{\grad f(z_b)\}$ \label{alg:add-row} 
			\ENDIF
			\IF{$f(z_b) < f(\hat{z})$}
			\STATE $\hat{z} \leftarrow z_b$ \label{alg:update}
			\ENDIF
	\ENDIF
 \ENDFOR
 \UNTIL{$|\Gamma|$ did not increase.}
 \RETURN{$\hat{z}$} \label{alg:return}
\caption{The algorithm}
\label{alg:algorithm}
\end{algorithmic}
\end{algorithm}

\bigskip

Before we show correctness, let us derive a bound on $\norm{H}_\infty$.
\begin{lem}
\label{lem:small-matrix}
In the course of the algorithm we have 
\[
\norm{H}_\infty \leq 2 \delta L.
\]
\end{lem}
\begin{proof}
By Lemma~\ref{lem:small-grads}, Line~\ref{alg:restrict-to-face} in the algorithm ensures $\norm{\grad f(x^\star)}_{\infty} \leq \max \{(\delta - 1) L, 0\}$, and we have $\norm{H}_\infty \leq \delta L$ when reaching Line~\ref{alg:end-init} for the first time.
By Line~\ref{alg:ip}, we only add a gradient $h = \grad f(z)$ to $\Gamma$ if $z$ satisfies $\norm{z - x^\star}_1 \leq \delta$.
Since the gradient of $g$ is Lipschitz-continuous the estimate $\norm{h}_\infty \leq \norm{\grad f(x^\star)}_{\infty} + \delta L \leq 2 \delta L$ holds and thus $\norm{H}_\infty \leq 2 \delta L$ in the course of the algorithm.
\end{proof}
We are ready to prove correctness of the algorithm.
\begin{lem}
\label{lem:correctness}
The algorithm solves at most $m (4 \delta^2 L + 1)^m$ integer programs.
The output is correct.
\end{lem}
\begin{proof}
By Line~\ref{alg:restrict-to-face} and Lemma~\ref{lem:small-grads} we can assume $\| \grad f (x^\star) \|_\infty \leq \max \{ (\delta - 1) L, 0 \}$ throughout the proof.
Since $\grad f(x) = (\grad g(Wx))^\intercal W$ every gradient we consider is a linear combination of the rows of $W$.
As we only add linearly independent vectors to $\Gamma$ we have $|\Gamma| \leq \rk (W) \leq m$.
Therefore the loop in Line~\ref{alg:outerloop} is repeated at most $m$ times.
Within the loop of Line~\ref{alg:b-loop}, at most $(2 \norm{H}_\infty \delta + 1)^{|\Gamma|}$ vectors $b$ are considered.
From Lemma~\ref{lem:small-matrix}, it follows that in total the algorithm solves at most $m (4 \delta^2 L + 1)^{m}$ integer programs in Line~\ref{alg:ip} before it terminates.

In order to show that the solution returned in Line~\ref{alg:return} is minimal, let $z^\star$ be a minimal integer solution with $\norm{z^\star - x^\star}_1 \leq \delta$.
This implies that the right-hand side $\hat{b} \df H z^\star$ will be considered in the loop of Line~\ref{alg:b-loop}.

When solving the system $H z = \hat{b}$ in Line~\ref{alg:ip}, assume the algorithm outputs another point $\hat{z} \in \{0,1\}^n$, $\hat{z} \neq z^\star$.
Since $\hat{z}$ minimizes the $\ell_1$-distance to $x^\star$ and $\norm{z^\star - x^\star}_1 \leq \delta$, we have $\norm{\hat{z} - x^\star}_1 \leq \delta$ in Line~\ref{alg:ip}.
If $\grad f (\hat{z})$ is linearly independent of $\Gamma$, we add it to $\Gamma$ and repeat eventually from Line~\ref{alg:outerloop}.
After at most $m$ iterations $\grad f (\hat{z})$ must depend linearly on $\Gamma$.
As $H$ comprises the vectors in $\Gamma$ as rows, there is a vector $y \in \R^r$ such that
$\grad f(\hat{z})^\intercal = y^\intercal H$.
This implies that
$$
\grad f(\hat{z})^\intercal (z^\star - \hat{z}) = y^\intercal H (z^\star - \hat{z}) = 0.
$$
Hence $f(z^\star) \geq f(\hat{z}) + \grad f(\hat{z})^\intercal (z^\star - \hat{z}) = f(\hat{z})$,
i.e., the solution $\hat{z}$ returned in Line~\ref{alg:return} is at least as good as $z^\star$, and therefore minimal.
\end{proof}

The previous arguments allow us to show our main result.

\begin{thm}[Convex minimization via solving IPs]
\label{thm:main}
Let $W \in \Z^{m \times n}$,
and $g$ be a convex function with $L$-Lipschitz continuous gradients such that $f(x) = g(Wx)$ admits a $\delta$-proximity.

Given a fractional minimum $x^\star \in [0,1]^n$ with at most $m$ fractional entries, $\delta$, and zeroth and first order oracle access to $f$, 
we can find an integer minimum $z^\star \in \{0,1\}^n$ by solving at most $m (4 \delta^2 L + 1)^m$ integer programs with $n$ binary variables, at most $m$ equality constraints, and the coefficients in the constraint matrix bounded by $2 \delta L$.
\end{thm}

Compared to the informal statement in Theorem~\ref{thm:informal},
we still require that $x^\star$ has at most $m$ fractional entries.
In the discussion after Lemma~\ref{lem:frac-entries} we saw that any constrained continuous minimum $x^\star$ suffices, provided we know $W$.
If we do not know $W$, we can still find such a minimum, provided we can find a continuous constrained minimum of a convex function.
\begin{lem}[Construction of $x^\star$ with few fractional entries]
\label{lem:frac-opt}
Assume we have a subroutine $\mathcal{A}$ that for any face $F \subseteq [0,1]^n$ computes a continuous constrained minimum $\hat{x} \in \argmin \{g(Wx): \ x \in F\}$.
We can find a continuous constrained minimum $x^\star \in [0,1]^n$ with at most $m$ fractional entries using at most $2n$ calls to the subroutine $\mathcal{A}$.
\end{lem}
\begin{proof}
We start by setting $F^{(0)} = [0,1]^n$ and computing 
\[
f_{\min} \df \min \{f (x): \ x \in F^{(0)}\},
\]
using the provided subroutine.

For $k=1,\dots,n$, we do the following.
In the $k$-th iteration, we consider the two faces $F^{(k-1)}_0 \df \{x \in F^{(k-1)}: \ x_k = 0\} \subseteq F^{(k-1)}$ and $F^{(k-1)}_1 \df \{x \in F^{(k-1)}: \ x_k = 1\} \subseteq F^{(k-1)}$,
and use the subroutine for solving the two problems
\begin{align*}
x^0 &\in \argmin \left\{ f(x) : x \in F^{(k-1)}_0 \right\}, \\
x^1 &\in \argmin \left\{ f(x) : x \in F^{(k-1)}_1 \right\}. \\
\end{align*}
If $f(x^0) = f_{\min}$, then set $F^{(k)} \df F^{(k-1)}_0$.
Similarly, if $f(x^0) > f_{\min}$ but $f(x^1) = f_{\min}$, we define $F^{(k)} \df F^{(k-1)}_1$.
If both $f(x^0) > f_{\min}$ and $f(x^1) > f_{\min}$ we simply define $F^{(k)} \df F^{(k-1)}$.

We claim that $\dim(F^{(n)}) \leq m$.
Indeed, by construction there exists
\[
\hat{x} \in \argmin \left\{ f(x) : x \in F^{(n)} \right\} \quad \text{with} \quad f(\hat{x}) = f_{\min}.
\]
Let $v$ be a vertex of the polyhedron $P \df \{x \in F^{(n)}: \ Wx = W \hat{x} \}$.
Then $v$ satisfies $f(v) = g(Wv) = g(W \hat{x}) = f(\hat{x}) = f_{\min}$ and by basic LP theory, $v$ has at most $m$ fractional entries.
Let $k$ be an index with $v_k \in \{0,1\}$.
We have $v \in F^{(n)} \subseteq F^{(k-1)}$, and since $v_k$ is integral, $v \in F^{(k-1)}_0$ or $v \in F^{(k-1)}_1$,
implying that $\dim(F^{(k)}) = \dim(F^{(k-1)})-1$.
As $v$ has at least $n-m$ integral entries, $\dim(F^{(n)}) \leq m$.
Consequently, any optimal solution in $F^{(n)}$ has at most $m$ fractional entries.
\end{proof}

Before we continue with showing concrete bounds on $\delta$ for separable convex and sharp functions, let us again comment on a potential improvement on the exponent in the running time to $\mathcal{O}(m^2)$ when $W$ is presented explicitly.

\begin{rem}
\label{rem:exponent}
As before one enumerates all points $b \in \Z^m$ in the projected space $H \R^n$ that are close to the image $Hx^\star$ of the continuous minimum $x^\star$.
For each point $b$, an integer program tests whether there is a feasible preimage $z \in \{0,1\}^n$, i.e., $b = Hz$.
If we know the matrix $W$ we can work in the space $W \R^n$ right away and search for $z$ such that $Wz = b$.
As the entries in $W$ are bounded by $\|W\|_{\infty}$ in contrast to $\|H\|_{\infty} \lesssim (m\|W\|_\infty)^m L$, the integer program in Line~\ref{alg:ip} can be solved more efficiently, saving a factor of $m$ in the exponent.

A crucial argument why the algorithm works without explicit knowledge of $W$ is the fact that the rows of $H$ are linear combinations of the rows of $W$.
In particular, if $\rk (H) = \rk (W)$, we have $\ker(H) = \ker(W)$.
As a consequence, if there are several solutions $x_1,x_2$ to the integer program $\{x \in \{0,1\}^n: \ Hx = b\}$, then $x_2 - x_1 \in \ker(W)$ and hence $g(Wx_1) = g(Wx_2)$, i.e., it does not matter which point $x_1$ or $x_2$ is returned in Line~\ref{alg:ip}.
In general, we only know $\ker(H) \supseteq \ker(W)$, but if the preimage returned by the integer program is not a global minimizer, then we can increase the rank of $H$.
\end{rem}

\section{Proximity results for classes of convex functions}
\label{sec:proximity}

This section is devoted to an analysis of upper bounds on $\delta$ for separable convex and sharp convex functions.
This will then allow us to make the running time of Algorithm~\ref{alg:algorithm} explicit.

\subsection{Separable convex functions}
\label{sec:proximity-sep-conv}
In this section, we assume in addition to the properties in Section~\ref{sec:intro} that $g$ is separable convex, i.e.,
$g(v) = \sum_{i=1}^m g_i(v_i)$ and each $g_i : \R \rightarrow \R$ is convex.
Our target is to derive an upper bound on $\delta$, see Definition~\ref{def:proximity}.
As it will turn out we will be able to derive a bound on $\delta$ that is independent on the dimension $n$ and improves a bound of~\citet{hs}.

When considering standard form integer linear optimization,~\citet{ew} show that proximity is also independent of the dimension.
In contrast to~\citet{hs}, not only the $\ell_{\infty}$ distance, but even the $\ell_1$-distance can be bounded.
As the technique of~\citet{ew} extends to separable convex functions, we can use it in combination with the special structure of our problem.

\begin{thm}[Proximity for separable convex functions]
\label{thm:proximity-sep-conv}
Let $g: \R^m \rightarrow \R$ be separable convex, i.e., $g(v) = \sum_{i=1}^m g_i(v_i)$ with each $g_i: \R \rightarrow \R$ convex, and $W\in \Z^{m \times n}$.
Then, $f(x) = g(Wx)$ admits a $\delta$-proximity with
\[
\delta \leq 2m ( 2 m \norm{W}_\infty + 1)^{m}.
\]
\end{thm}
\begin{proof}
Let $x^\star \in [0,1]^n$ be a constrained continuous minimum, and assume that $x^\star$ has at most $m$ fractional entries by Lemma~\ref{lem:frac-entries}.
Consider the integer program
\begin{alignat}{2}
\label{eq:ip-reformulation}
\min \ & & g(y) \\
s.t. \ & \quad & Wx &= y, \nonumber \\
& & x &\in \{0,1\}^n, \nonumber \\
& & -n \norm{W}_{\infty} \leq y_j &\leq n \norm{W}_{\infty}, \qquad j=1,\dots,m, \nonumber \\
& & y &\in \Z^m, \nonumber 
\end{alignat}
and observe that it is equivalent to our initial problem $\min \{ f(x): \ x \in \{0,1\}^n\}$.
Thus, $\binom{x^\star}{Wx^\star}$ is a minimal solution to the relaxation with at most $2m$ fractional entries.
An extension of the proximity result of~\citet{ew} to a separable convex function $h: \R^n \rightarrow \R$ can be found for instance in~\citet[Prop.~$2.24$]{thesis}, stating the following.

Let $A \in \Z^{m \times n}$ and $x^\star$ be a minimal solution of 
$\min \{h(x):\ Ax = b, \ \ell \leq x \leq u\}$
that has $\varphi \in [n]$ fractional entries.
If the problem
$\min \{h(x):\ Ax = b, \ \ell \leq x \leq u, \ x \in \Z^n\}$
is feasible, then there exists a minimal solution $z^\star$ satisfying
$\|x^\star - z^\star\|_1 \leq \varphi \gamma(A)$,
where $\gamma(A)$ is an upper bound on the length of the \emph{Graver basis elements} of the matrix $A$.
Lemma $2.4$ of \citet{thesis} states that $\gamma(A) \leq (2m\|A\|_{\infty}+1)^m$.

Applied to our setting,
there exists a minimal integral solution $\binom{z^\star}{Wz^\star}$ to~\eqref{eq:ip-reformulation} satisfying
$
\norm{\binom{x^\star}{Wx^\star} - \binom{z^\star}{Wz^\star}}_1 \leq 2m(2m\norm{W}_\infty + 1)^m
$.
Since $\norm{\binom{x^\star}{Wx^\star} - \binom{z^\star}{Wz^\star}}_1 = \norm{x^\star - z^\star}_1 + \norm{W(x^\star - z^\star)}_1$, the result follows.
\end{proof}

In terms of running time, we can combine this result with Theorem~\ref{thm:main} and apply
the algorithm of~\citet{ew}.

\begin{cor}
\label{cor:sep-conv}
Let $W \in \Z^{m \times n}$,
and $g$ be a separable convex function with $L$-Lipschitz continuous gradients, such that $\grad g(v) \in \Z^{m}$ whenever $v \in \Z^m$.

Given $m$, $\norm{W}_{\infty}$, $L$, and a constrained continuous minimum $x^\star \in [0,1]^n$ with at most $m$ fractional entries,
there exists an algorithm computing an integer minimum $z^\star \in \{0,1\}^n$ for $f$
in $\poly(n) (Lm \norm{W}_{\infty} + 1)^{ \mathcal{O} (m^3) }$ arithmetic operations and $(Lm \norm{W}_{\infty} + 1)^{ \mathcal{O} (m^2) }$ first order oracle calls.
\end{cor}

\subsection{Sharp convex functions}
\label{sec:proximity-sharp}
In this section, let $g$ still satisfy the assumptions made in Section~\ref{sec:intro}.
In addition, we assume that $g$ satisfies the \emph{H\"older error bound (HEB) condition},
which is also referred to as \emph{sharpness condition}.
This is a generalization of strong convexity; whereas strong convexity is a global property, sharpness is local, i.e., only a bound towards the minimum is ensured, but not between any two points.
We point out that there are varying definitions of sharpness in the literature, and settle for the following.
\begin{definition}[Sharpness]
\label{def:shapness}
Let $P \subseteq \R^m$ be a convex set, $g: P \rightarrow \R$ be convex, and assume $v^\star \in P$ with $g(v^\star) = \min_{v \in P} g(v)$ exists.

The function $g$ satisfies the \emph{H\"older error bound condition} (also: the \emph{sharpness condition}) on $P$ with parameters $0 < \mu < \infty$ and $\theta \in [0,1]$, if the following holds for all $u \in P$:
\begin{align}
\label{eq:sharpness}
\mu \left( g(u) - g(v^\star) \right)^\theta \geq \norm{u-v^\star}_1.
\end{align}
\end{definition}
Since we are interested in minimizing $f(x) = g(Wx)$ with $x \in \{0,1\}^n$ and $W \in \Z^{m \times n}$, we set $P = \{Wx : \ x \in [0,1]^n\} \subseteq [-n\norm{W}_\infty, n\norm{W}_\infty]^m$.
Though this definition implies uniqueness for a continuous constrained minimum of $g$, the function $f(x) = g(Wx)$ may have several minima in general.

We are now prepared to show proximity for sharp convex functions.

\begin{thm}[Proximity for sharp convex functions]
\label{thm:proximity-sharp}
Let $W \in \Z^{m \times n}$, $g: \{Wx : \ x \in [0,1]^n \} \rightarrow \R$ be convex, have $L$-Lipschitz continuous gradients, and satisfy the sharpness condition with parameters $\mu>0$ and $\theta \in [0,1]$.
Then $f(x) = g(Wx)$ admits a $\delta$-proximity with
\[
\delta \leq \left( \mu \left( m^4 \tfrac{L}{4} \norm{W}_\infty^2 \right)^\theta + m + 2 \right) (2m \norm{W}_\infty + 1)^m.
\]
\end{thm}
\begin{proof}
Let $x^\star \in [0,1]^n$ be a constrained continuous minimal solution for $f(x) = g(Wx)$ over $[0,1]^n$ with at most $m$ fractional entries by Lemma~\ref{lem:frac-entries}, and $z^\star \in \{0,1\}^n$ a minimal integral solution that is closest to $x^\star$.
We will first estimate the distance between $v \df Wx^\star$ and $u \df Wz^\star$.
In combination with the sensitivity result (Lemma~\ref{lem:sensitivity}), the bound for $x^\star$ and $z^\star$ will follow.
To this end, let $w \df W\lfloor x^\star \rceil$, where $\lfloor x^\star \rceil$ denotes $x^\star$ rounded component-wise to the nearest integer, and observe that $v$ is minimal for $\min\{g(y): \ y \in W[0,1]^n\}$.

If $x^\star_k$ is fractional, both $v \pm \varepsilon W_{\cdot k}$ are feasible.
Lemma~\ref{lem:first-order-opt} implies that 
$$
\grad g(v)^\intercal W_{\cdot k} = 0,
$$
and thus $\grad g(v)^\intercal (w - v) = 0$.
Using the Lipschitz-continuity of the gradients, we have
\begin{align*}
g(w) - g(v) &\leq \grad g(w)^\intercal (w-v) \\
	&= \grad g (v)^\intercal (w-v) - (\grad g (v) - \grad g (w))^\intercal (w-v) \\
	&\leq \norm{\grad g(v) - \grad g(w)}_\infty \norm{w-v}_1 \\
	&\leq L \norm{w-v}^2_1 \\
%g(w) - g(v) &\leq \grad g(v)^\intercal (w - v) + L \norm{w-v}_1^2 \\
	&= L \norm{W(\lceil x^\star \rfloor - x^\star)}_1^2 \nonumber \\
	&\leq L (m \norm{W ( \lceil x^\star \rfloor - x^\star)}_{\infty})^2 \nonumber \\
	&\leq L (m \norm{W}_{\infty} \norm{ \lceil x^\star \rfloor - x^\star }_{1})^2 \nonumber \\
	&\leq m^4 \norm{W}_\infty^2 \tfrac{L}{4}. \nonumber 
\end{align*}
On the other hand, sharpness and uniqueness of $v$ allow us to estimate
\[
g(u) - g(v) \geq \left( \tfrac{1}{\mu} \norm{u-v}_1 \right)^{1/\theta}.
\]
Combining the two established inequalities, we obtain
\[
m^4 \norm{W}_\infty^2 \tfrac{L}{4} \geq g(w) - g(v) \geq g(u) - g(v) \geq \left( \tfrac{1}{\mu} \norm{u-v}_1 \right)^{1/\theta},
\]
and rearranging yields $\norm{u-v}_1 \leq \mu \left( m^4 \tfrac{L}{4} \norm{W}_\infty^2 \right)^\theta$.
Now, by Lemma~\ref{lem:sensitivity}, we see that $z^\star$ satisfies
\[
\norm{x^\star - z^\star}_1 \leq \left( \mu \left( m^4 \tfrac{L}{4} \norm{W}_\infty^2 \right)^\theta + m + 2 \right) (2m \norm{W}_\infty + 1)^m.
\qedhere
\]
\end{proof}

A combination of this result with Theorem~\ref{thm:main} and
the algorithm of~\citet{ew} gives the following running time.

\begin{cor}
\label{cor:sharp-conv}
Let $f(x) = g(Wx)$ for $W \in \Z^{m \times n}$,
and $g$ be a convex function that has $L$-Lipschitz continuous gradients, 
such that $\grad g(v) \in \Z^{m}$ whenever $v \in \Z^m$.
Moreover, let $g$ be sharp with parameters $\mu$ and $\theta$.
%, and have a unique continuous constrained minimum on $\{Wx: \ x \in [0,1]^n\}$.

Given $m$,$\norm{W}_{\infty}$,$L$,$\mu$,$\theta$ and a constrained continuous minimum $x^\star \in [0,1]^n$ with at most $m$ fractional entries, there exists an algorithm computing an integer minimum $z^\star \in \{0,1\}^n$ for $f$
in 
\[
\poly(n) \left( \mu \left( m^4 L \norm{W}_\infty^2 \right)^\theta + 1\right)^{\mathcal{O}(m^3)}
\]
arithmetic operations and $\left( \mu \left( m^4 L \norm{W}_\infty^2 \right)^\theta + 1\right)^{\mathcal{O}(m^2)}$ first order oracle calls.
\end{cor}

% ---------------------------------bibliography-----------------------------
\bibliographystyle{abbrvnat}
\bibliography{bibliography}

\end{document}